\documentclass[a4paper,12pt,twoside]{article}
\usepackage{amsthm,amssymb,amsmath,amscd, verbatim,url,amsfonts}
\newtheorem{thm}{Theorem}[section]

\theoremstyle{definition}

\def\li{\textrm{li}}

\begin{document}
\title{A new bound for the smallest $x$ with $\pi(x) > \li(x)$}
\author{Kuok Fai Chao and Roger Plymen}
\date{}
\maketitle \noindent School of Mathematics, Alan Turing building,
Manchester University\\ Manchester M13 9PL, England\\

kchao@maths.manchester.ac.uk, plymen@manchester.ac.uk\\

\begin{abstract}
 We reduce the leading term in Lehman's theorem.
This improved estimate allows us to refine the main theorem of
Bays \& Hudson \cite{BH}. Entering $2,000,000$ Riemann zeros, we
prove that there exists $x$ in the interval
$[\textmd{exp(727.951858)}, \textmd{exp(727.952178)}]$ for which
$\pi(x)-\li(x) > 3.2 \times 10^{151}$.  There are at least
$10^{154}$ successive integers $x$ in this interval for which
$\pi(x)>\li(x)$. This interval is strictly a sub-interval of the
interval in Bays \& Hudson, and is narrower by a factor of about
$12$.
\end{abstract}



\section{Introduction}
Let $\pi(x)$ denote the number of primes less than or equal to
$x$, and let $\li(x)$ denote the logarithmic integral. The
notation $f(x) = \Omega_{\pm}g(x)$ means that \[\lim \sup_{x \to
\infty} f(x)/g(x) > 0, \quad \quad \lim \inf_{x \to \infty}
f(x)/g(x) <0
\] There was, in 1914, overwhelming numerical evidence that
$\pi(x)<\li(x)$ for all $x$. In spite of this, Littlewood
\cite{Li} announced that
\[\pi(x) - \li(x) = \Omega_{\pm}(x^{1/2}(\log x)^{-1} \log \log \log
x)\]
 This implies
that $\pi(x) - \li(x)$ changes sign infinitely often. Littlewood's
method provided, even in principle, no definite number $X$ before
which $\pi(x) - \li(x)$ changes sign. For a recent proof of
Littlewood's theorem, see \cite[Theorem 15.11]{M}.

In the course of the 20th century, successive numerical upper
bounds were found by Skewes \cite{S0}, Skewes \cite{S}, Lehman
\cite{L}, te Riele \cite{R}. For Littlewood's own account of the
discovery of the Skewes numbers, see \cite[p. 110--112]{L2}.

The smallest value of $x$ with $\pi(x) \geq \li(x)$ will be
denoted $\Xi$, as in the recent paper by Kotnik \cite{K}.  In the
course of a systematic computational study, Kotnik proves that
\[10^{14} < \Xi.\]

We now explain the main idea in \cite{L}.  Lehman's theorem is an
\emph{integrated version of the Riemann explicit formula}. His
method was to integrate the function $u \mapsto \pi(e^u) -
\li(e^u)$ against a Gaussian kernel over a carefully chosen
interval $[\omega - \eta, \omega + \eta]$. The definite integral
so obtained is denoted $I(\omega,\eta)$. Let $\rho = 1/2 +
i\gamma$ denote a Riemann zero with $\gamma > 0$ and let
\[H(T,\omega): = - 2 \Re \sum_{0 < \gamma \leq T}\frac{e^{i\gamma
\omega}}{\rho} e^{-\gamma^2/2\alpha}.\]  The $\alpha$ in this
formula is related to the kernel chosen. Lehman proved the
following equality
\[I(\omega,\eta) = -1 + H(T,\omega) + R
\] together
with an explicit estimate $|R| \leq \epsilon$.  This creates the
inequality \[ I(\omega,\eta) \geq H(T,\omega) - (1
+\epsilon).\]The problem now is to prove that
\begin{eqnarray}
\label{Sigma}  H(T,\omega) > 1 + \epsilon.
\end{eqnarray} If (\ref{Sigma}) holds, then
$I(\omega,\eta)>0$ and so there exists $x \in [e^{\omega - \eta},
e^{\omega + \eta}]$ for which $\pi(x) > \li(x)$.  In order to
establish (\ref{Sigma}), numerical values of the Riemann zeros
with $|\gamma| < T$ are required.  Each term in $H(T,\omega)$ is a
complex number determined by a Riemann zero.  It is necessary that
the real parts of these complex numbers, which are spiralling
towards $0$, reinforce each other sufficiently for (\ref{Sigma})
to hold. The only known way of establishing this is by numerical
computation.  When $T$ is large, this requires a computer.

In 2000, Bays \& Hudson \cite{BH} made the following selection:
\[
\omega = 727.95209, \quad \eta = 0.002.\]
 The
interval itself is $[\textmd{exp(727.95009)},
\textmd{exp(727.95409)}].$ We note that the numerical value of
$\exp(727.95409)$ is \emph{incorrectly stated} in \cite[Theorem
2]{BH}.

In this article, we reduce the leading term in Lehman's theorem.
This enables us to select the following parameters:
\[
\omega = 727.952018, \quad \eta = 0.00016.\]
 The
interval is $[\textmd{exp(727.951858)}, \textmd{exp(727.952178)}]$
and so an upper bound for the first crossover $\Xi$ is
\[\Xi < \textmd{exp(727.952178)} < 1.398344 \times
10^{316}\]

Our interval is strictly a sub-interval of the Bays-Hudson
interval. It is narrower by a factor of about $12$, and creates
the smallest known upper bound.

The function $H(T,\omega)$ is an initial part of the series
\[H(\omega): = - 2 \Re \sum_{0 < \gamma}\frac{e^{i\gamma \omega}}{\rho} e^{-\gamma^2/2\alpha}.\]

As Rademacher observed in 1956 \cite{Ra}, the Riemann Hypothesis
plus Weyl's criterion imply that, for each $\omega > 0$,  the
sequence
\[\{\exp(i \gamma\omega): \zeta(1/2 + i\gamma) = 0, \gamma
>0\}\] is \emph{equidistributed} in the unit circle. So we may
expect a fair amount of cancellation to take place in the series
$H(\omega)$.  This may help to explain why it is so difficult to
find a number $\omega$ for which $H(T,\omega)$ exceeds $1$.

We reflect, for a moment, on the Weil explicit formula. This is an
identity between two distributions \cite[p.39]{Pa}. It is well
established that certain classical explicit formulas follow from
the Weil explicit formula, by picking suitable test-functions. For
example, classical formulas for Dirichlet $L$-series may be
derived in this way, see \cite[Theorem 3.2, p.340]{Lang}. We are
led to ask whether the Lehman formula can be obtained from the
Weil explicit formula by picking a suitable test function. We hope
to pursue this idea elsewhere.

We would like to thank Andrew Odlyzko for supplying us with the
first $2,000,000$ Riemann zeros, Jon Keating for drawing our
attention to \cite{Ra},  Christine Lee for improvements in the
exposition, Aleksandar Ivi\'c for drawing our attention to the
inequalities of Panaitopol \cite{P}, and Nick Gresham for
assistance with numerical calculations. Finally, we thank the
referee for his many detailed and constructive comments.

\section{The Leading Term}

We begin this section with Lehman's theorem.

\begin{thm}\label{ingham} (Lehman \cite{L}) Let $A$ be a positive number such that
$\beta=\frac{1}{2}$ for all zeros $\rho=\beta+i\gamma$ of
$\zeta(s)$ for which $0<\gamma\leq A$. Let $\alpha$, $\eta$ and
$\omega$ be positive numbers such that $\omega-\eta>1$ and
\begin{equation}\label{c1}
2/A\leq 2A/\alpha\leq\eta\leq \omega/2.
\end{equation}
Let
\begin{equation*}\label{K}
K(y):=\sqrt{\frac{\alpha}{2\pi}}e^{-\alpha y^{2}/2}
\end{equation*}
\begin{equation*}\label{omega}
I(\omega,\eta):=\int_{\omega-\eta}^{\omega+\eta}K(u-\omega)ue^{-u/2}\{\pi(e^{u})-\textrm{li}(e^{u})\}du
\end{equation*}
Then for $2\pi e<T\leq A$
\begin{equation*}\label{I}
I(\omega, \eta) = -1-\sum_{0<|\gamma|\leq
T}\frac{e^{i\gamma\omega}}{\rho}e^{-\gamma^{2}/2\alpha}+R
\end{equation*}
where
\begin{equation*}
|R|\leq  s_1+s_2+s_3+s_4+s_5+s_6
\end{equation*}
with
\begin{eqnarray*}
s_{1}&=&\frac{3.05}{\omega-\eta}\\
s_{2}&=&4(\omega+\eta)e^{-(\omega-\eta)/6}\\
s_{3}&=&\frac{2e^{-\alpha\eta^{2}/2}}{\sqrt{2\pi \alpha }\eta}\\
s_{4}&=&0.08\sqrt{\alpha}e^{-\alpha\eta^{2}/2}\\
s_{5}&=&e^{-T^{2}/2\alpha}\{\frac{\alpha}{\pi
T^{2}}\log\frac{T}{2\pi} + \frac{8\log T}{T} +
\frac{4\alpha}{T^{3}}\}\\
s_{6}&=&A\log Ae^{-A^{2}/2\alpha+(\omega + \eta)/2}
\{4\alpha^{-1/2}+15\eta\}
\end{eqnarray*}
If the Riemann Hypothesis holds, then conditions (\ref{c1}) and
the term $s_6$ in the estimate for $R$ may be omitted.
\end{thm}
\medskip

 For the rest of the paper $\rho = \beta + i\gamma$
will denote a zero of the Riemann zeta function $\zeta(s)$ for
which $0 < \beta < 1$. We will refine a part of Lehman's proof.
This allows us to reduce the term $s_{1}$ in Lehman's theorem.

The logarithmic integral is defined as follows \cite[p.82]{I}:
\[
\li(e^z): = \int_{- \infty + iy}^{x + iy}\frac{e^t}{t}dt\] where
$z = x+iy, \, y \neq 0$.  For $x
> 1$,\, $\li(x)$ is then defined as follows:
\[
\li(x): = \frac{1}{2}[\li(x+i0) + \li(x - i0)]
\] In this way, we recover the classical definition of $\li(x)$ as
an integral principal value \cite[p.82]{I}:
\[
\li(x) = \lim_{\epsilon \to 0}\left(\int_0^{1-\epsilon}
\frac{e^t}{t}dt + \int_{1+\epsilon}^x \frac{e^t}{t}dt\right)\] For
a detailed account of the logarithmic integral, see
\cite[p.38--41]{Leb}.

 We define
\begin{equation*}
J(x): = \pi(x) + \frac{1}{2} \pi(x^{1/2}) + \frac{1}{3}
\pi(x^{1/3}) \ldots \, \label{Lex1}
\end{equation*}
and recall the Riemann-von Mangoldt explicit formula:
\begin{equation}
J(x) = \li(x) - \sum_{\rho} \li(x^{\rho}) + \int_{x}^{\infty}
\frac{du}{(u^2 - 1) u \log u} - \log 2\ \label{Lex2}
\end{equation}
valid for $x > 1$. According to \cite[(3.2) and (3.6), p.69]{RS}
we have, for all $x > 1$,
\begin{align}
&\pi(x) = \frac{x}{\log x} + \frac{3\vartheta_1(x)x/2}{(\log x)^2} \notag \\
&\pi(x) = \frac{\vartheta_2(x) 2x}{\log x}  \label{Lpi}
\end{align}
\label{varthetas}with $|\vartheta_1(x)| < 1,\,\vartheta_2(x) <
0.62753$. There are at most $[(\log x) /\log 2]$ terms in $J(x)$.
This allows us implicitly to define $\vartheta_3(x)$ by the
following equation: \begin{align} J(x) = \pi(x) + \frac{1}{2}
\pi(x^{1/2}) + \frac{1}{3} \vartheta_3(x) \pi(x^{1/3})
\left(\frac{\log x}{\log 2} \right) \end{align} Then we have
$\vartheta_3(x) < 1$. Combining (4) and (6), we have
\begin{align*}
\pi(x) - \li(x) &= - \sum_{\rho} \li(x^{\rho}) + \int_{x}^{\infty} \frac{du}{(u^2 - 1) u \log u} - \log 2 \\
                   & - \frac{1}{2} \pi(x^{1/2}) -  \frac{1}{3} \vartheta_3(x) \pi(x^{1/3}) \left(\frac{\log x}{\log 2} \right)\,.
\end{align*}
Substituting the first expression in \eqref{Lpi} for
$\pi(x^{1/2})$ and the second expression for $\pi(x^{1/3})$, we
get \begin{align*}
\pi(x) - \li(x) &= - \sum_{\rho} \li(x^{\rho}) + \int_{x}^{\infty} \frac{du}{(u^2 - 1) u \log u} - \log 2 \\
                   & - \left(\frac{x^{1/2}}{\log x} + 3\frac{\vartheta_1(x^{1/2})x^{1/2}}{(\log x)^2} \right) -
                   \vartheta_3(x) \left(\frac{\vartheta_2(x^{1/3})2 x^{1/3}}{\log 2} \right).
\end{align*} Let \[\vartheta_4(x): = \int_{x}^{\infty} \frac{du}{(u^2-1)u \log u} - \log
2\] For $x > 2$ we have the following bounds:
\[- \log 2 < \vartheta_4(x)  < 1/2 - \log 2.
\]
We now have
\begin{align*}
\pi(x) - \li(x) &= - \sum_{\rho} \li(x^{\rho})  -
\frac{x^{1/2}}{\log x} \\ & +
3\frac{\vartheta_1(x^{1/2})x^{1/2}}{(\log x)^2} + \vartheta_4(x) -
                   \frac{\vartheta_3(x)\vartheta_2(x^{1/3})2 x^{1/3}}{\log 2}
\end{align*}Now define $\vartheta(x)$ as follows:
\[
4\vartheta(x)x^{1/3}: = \vartheta_4(x)  -
\frac{\vartheta_3(x)\vartheta_2(x^{1/3})2x^{1/3}}{\log 2}\] Then
we have
\[
|\vartheta(x)| < 4^{-1}x^{-1/3}(1/2 - \log 2) + 1/2 < 1 - \log 2 <
1\] for all $x > 2$.
 Here is where our method differs from
Lehman's approach: we keep the estimate $\vartheta_1(x)$ separate
from $\vartheta(x)$. We have
\begin{equation}\label{main}
\pi(x) - \li(x) = -\sum_{\rho} \li(x^{\rho}) - \frac{x^{1/2}}{\log
x} + \frac{3\vartheta_1(x^{1/2})x^{1/2}}{(\log x)^2} +
\vartheta(x) 4x^{1/3}
\end{equation}

\noindent Now we improve the bound for $\vartheta_1(x)$. We quote
a result of Panaitopol \cite[Theorem 1]{P}:
\begin{equation}\label{kkk}
\pi(x)<\frac{x}{\log x-1-(\log x)^{-0.5}} \quad \textrm{for all}
\quad x\geq 6
\end{equation}
From (\ref{Lpi}) and (\ref{kkk}), we get
\begin{equation}\label{theta}
\frac{x}{\log x}+\frac{\frac{3}{2}
x\vartheta_1(x)}{\log^{2}x}<\frac{x}{\log x-1-(\log x)^{-0.5}}
\end{equation}
Denote $y= y(x): = (\log x)^{\frac{1}{2}}$. The inequality
(\ref{theta}) will lead to an upper bound for $\vartheta_1(x)$:
\begin{eqnarray}\label{yyy}
0<\vartheta_1(x) <\frac{2}{3}\cdot \frac{y^{3}+y^{2}}{y^{3}-y-1}
\quad \textrm{for all } x\geq 6.
\end{eqnarray}
We define \[F(y):=\frac{y^{3}+y^{2}}{y^{3}-y-1}.\] We have $F(y)
> 1,\;F(y) \to 1$ as $y \to \infty$, and
\[
F'(y) = (-4y^3 - 3y^2 -1)/(y^3 - y - 1)^2 < 0\] so that $F$ is a
monotone decreasing function.  By (\ref{yyy}) we have
\[\vartheta_1(e^v) < \frac{2}{3}\cdot F(\sqrt v)\]  and so
\[
3\vartheta_1(e^{u/2}) < 2 F(\sqrt{\tfrac{727}{2}}) < 2.1111\] if
$u \geq 727$.

 From (\ref{main}) we have immediately:
\begin{equation*}
ue^{-u/2}\{\pi(e^{u})-\textrm{li}(e^{u})\}=-1-\sum_{\rho}ue^{-u/2}\textrm{li}(e^{\rho
u})+\vartheta_1(e^{u/2})\frac{3}{u}+4\vartheta(e^u) ue^{-u/6}
\end{equation*}

Now $K$ is a standard Gaussian distribution, and so

\[\int_{\omega-\eta}^{\omega+\eta}K(u-\omega)du =
\int_{-\eta}^{\eta}K(v)dv < 1.
\]
If $\omega - \eta > 727$ then we have the estimate
\[
\left|\int_{\omega-\eta}^{\omega+\eta}K(u-\omega)(\vartheta_1(e^{u/2})
\frac{3}{u}+4\vartheta ue^{-u/6})du \right| \leq
\frac{2.1111}{\omega-\eta}+4(\omega+\eta)e^{-(\omega-\eta)/6}
\]
We have replaced the term $s_1$ by $s'_1$:
\begin{equation*}
s'_1=\frac{2.1111}{\omega-\eta}
\end{equation*}
Following the steps in Lehman's proof \cite{L}, we are led to a
new estimate $|R'|$:
\begin{equation}\label{new r}
|R'|\leq s'_1+s_2+s_3+s_4+s_5+s_6
\end{equation}

\begin{thm}\label{coffee}
Let $A$ be a positive number such that $\beta=\frac{1}{2}$ for all
zeros $\rho=\beta+i\gamma$ of $\zeta(s)$ for which $0<\gamma\leq
A$. Let $\alpha$, $\eta$ and $\omega$ be positive numbers such
that $\omega-\eta > 727$ and
\begin{equation}\label{con1}
2/A\leq 2A/\alpha\leq\eta\leq \omega/2.
\end{equation}
Let $K(y)$ and $I(\omega,\eta)$ be defined as in Theorem 2.1. Then
for $2\pi e<T\leq A$ we have
\[I(\omega,\eta) = -1-\sum_{0<|\gamma|<T}\frac{e^{i\gamma\omega}}{\rho}e^{-\gamma^{2}/2\alpha}+R'
\]
where an upper bound for $R'$ is given by (\ref{new r}).

 If the Riemann Hypothesis
holds, then conditions (\ref{con1}) and the term $s_6$ in the
estimate for $R'$ may be omitted.
\end{thm}\medskip

{\sc Note}. We are setting out to narrow the interval $[\omega -
\eta, \omega + \eta]$.  However, the inequalities (\ref{con1})
include a \emph{lower bound} on $\eta$.  We explain briefly how
this lower bound arises. Without the Riemann Hypothesis, Lehman
proves, by means of several intricate estimates
\cite[p.404--406]{L}, that the inequalities (\ref{con1}) are a
sufficient condition for the following crucial estimate:
\[
\left|\sum_{|\gamma| > A} \int_{\omega - \eta}^{\omega + \eta}
K(u-\omega) ue^{-u/2}\li(e^{\rho u})du\right| \leq s_6.\]\medskip

\section{Numerical Results}
In this section we exploit the reduced term $s'_{1}$ to obtain
improved numerical results. We commence with some remarks
concerning the accuracy of computing by machine. Adopting notation
similar to te Riele \cite{R} we set
\[
H(T, \alpha, \omega) = - \sum_{0 < |\gamma|\leq
T}e^{-\gamma^2/2\alpha}\cdot\frac{e^{i\omega\gamma}}{\rho}
     = - \sum_{0<\gamma\leq T}t(\gamma,\alpha, \omega)
\]
where
\[
 t(\gamma,\alpha,\omega) = e^{-\gamma^2/2\alpha}\cdot\frac{\cos(\omega\gamma)+2\gamma\sin(\omega\gamma)}{\tfrac{1}{4} + \gamma^2}
\]
Fixing $N = 2,000,000$, we denote by $\gamma^{*}_{i}$ the
approximations to the true values of $\gamma_{i}$, correct to nine
decimal places, computed by Odlyzko for $1\leq i\leq N$, and set
$T = 1131944.47182487 > \gamma_{N}$. Let $H^{*}(T, \alpha,\omega)$
be the value obtained by taking the sum up to $T$ using
$t(\gamma^{*},\alpha,\omega)$ and let $H^{*}_{M}(T, \alpha,
\omega)$ be the result of computing the same sum by machine. The
function $H^{*}_{M}$ therefore depends on many machine-specific
details, including the implementation of standard library
functions \texttt{sin}, \texttt{cos} and \texttt{exp}; similar
considerations apply to the evaluation of $s'_{1}, \ldots, s_{6}$.
In the case of $H$ the quantity we wish to control is
\begin{displaymath}
|H- H^{*}_{M}| \leq |H - H^{*}| + |H^{*}-H^{*}_{M}|
\end{displaymath}
The terms on the right represent errors arising from the
conditioning of $H$ upon the $\gamma_{i}$ and numerical
instability in computing by machine respectively.

Following \cite{R} we have for $\gamma < \alpha$
\begin{displaymath}
 |H(T, \alpha, \omega) - H^{*}(T, \alpha, \omega)| < \sum_{0<\gamma\leq T}|\gamma-\gamma^{*}|\cdot M(\gamma,\alpha,\omega)
\end{displaymath}
with
\begin{displaymath}
M(\gamma,\alpha,\omega) =
e^{-\gamma^{2}/2\alpha}\left(\frac{2\omega}{\gamma} +
\frac{\omega}{\gamma^2}+\frac{2}{\alpha}+\frac{2}{\gamma^3}+\frac{4}{\gamma^2}\right)<\frac{3\omega+8}{\gamma}
\end{displaymath}

Now $14 < \gamma < 1131945 < \alpha$, and in the region of
interest, we have $\omega < 728$, so we may use the estimate
\begin{displaymath}
 |H(T, \alpha, \omega) - H^{*}(T, \alpha, \omega)| < 10^{-9}\cdot2192\cdot\sum^{N}_{i=1}\frac{1}{\gamma_{i}}
\end{displaymath}
and we find numerically that $\sum^{N}_{i=1}\frac{1}{\gamma_{i}} <
12$, so that $|H-H^{*}|$ is bounded above by $3\times 10^{-5}$.

It remains to address the extent of the additional error arising
from machine computation of $H^{*}$ and the quantities
$s'_{1},\ldots,s_{6}$. Let us denote by $R'_{M}$ the
machine-evaluated sum $s'_{1}+\cdots+s_{6}$. Our initial
experiments were conducted using \textsc{Matlab}, but to speed
sampling of the space of parameters ($\omega, \alpha$, $\eta$,
$A$) we re-implemented matters in \texttt{C} on an \texttt{x86\_64
GNU/Linux} system using native double precision and routines from
the \texttt{GNU} standard mathematics library \texttt{libm}. The
results for specific valid choices of  $(\omega, \alpha, \eta, A)$
with $H^{*}_{M}-(1+R'_{M}) > 1\times10^{-4}$ were then re-computed
using the arbitrary-precision libraries
\textsc{arprec}\cite{ARPREC} and \textsc{mpfr}\cite{MPFR}, running
at up to 100 digits and 1024 bits of precision respectively. Upon
rounding to 7 decimal places all values obtained were in
agreement; we are therefore confident that the cumulative effects
of adverse numerical phenomena lie well below the threshold of
$1\times 10^{-6}$.

A simple strategy for selecting suitable $(\omega, \alpha, \eta,
A)$ is to make order of magnitude estimates of the exponential
factors in the terms $s_2,\ldots,s_6$. To start note that with two
million zeros and $\omega$ near $728$ the exponential factor in
$s_{5}$ will be of the same order of magnitude as the leading term
$s'_{1}$ when $\alpha$ is near $10^{11}$. A larger $\alpha$ means
that $\eta$ may be taken smaller, resulting in a narrower interval
$(\omega-\eta, \omega+\eta)$, but pushing $\alpha$ to $1.35\times
10^{11}$ seems to be as far as one may safely travel in this
direction. By repeated subdivision and scanning over subintervals
of $(727,728)$ with this $\alpha$ we find $H^{*}_{M} - 1$ can be
made to exceed the leading term $s'_1
> 2\times 10^{-3}$ near $\omega=727.95202$. To control the
contribution from $s_{6}$ we need to take an $A$ the same order of
magnitude as $(\alpha\omega)^{1/2}$, whereupon the constraint
$\alpha \leq A^2$ forces $A$ near $1\times10^{7}$. We cannot take
$\eta$ an order of magnitude smaller than $10^{-5}$ without losing
control of $s_{4}$, and the constraint $\alpha\eta\geq 2A$ forces
$\eta$ close to $2\times 10^{-4}$. To keep $H^{*}_{M} - 1 -
R'_{M}$ safely above $1\times 10^{-4}$ we chose the following
values:
\begin{align*}
\omega &= 727.952018 & \eta &= 0.00016 \\
 \alpha &= 1.34\times10^{11} &A &= 1.022\times 10^{7}
\end{align*}
Then we obtain $H^{*}_{M} > 1.006569$ and the estimates
\begin{align*}
s'_{1} &< 0.002901 & s_{2} &< 10^{-49}\\
s_{3} &< 10^{-746} & s_{4} &< 10^{-740} \\
 s_{5} &< 0.003380 & s_{6} &< 10^{-5}
\end{align*}
so that
\begin{eqnarray*}
\label{II} I(\omega,\eta) & \geq & H - (1+|R'|) \\
& \geq & H^{*}_{M} - (1+R'_{M}) -  3 \times 10^{-5} - 1 \times 10^{-6} \\
& \geq & 2\times 10^{-4} >0
\end{eqnarray*}
Thus there is a value of $u$ in the interval
\begin{displaymath}
(\omega-\eta,\omega+\eta) = (727.951858,727.952178)
\end{displaymath}
for which $\pi(e^{u})> \li(e^{u})$. Let
\[
 F(u) = ue^{-u/2}\cdot(\pi-\li)(e^u)
\]
then we have shown that
\[
I(\omega,\eta) = \int_{\omega-\eta}^{\omega+\eta}K(u-\omega)\cdot F(u)du
\geqslant\delta
\]
with $\delta = 2\times 10^{-4}$,  therefore since $K(u)du$ is a probability measure,
\[
0<\delta\leq \int_{\omega-\eta}^{\omega+\eta}K(u-\omega)\cdot\sup
F(u)du < \sup F(u)
\]
Now $F(u)$ is continuous except where $u$ happens to be the
logarithm of a prime, so it follows that for $u$ in some
subinterval of $(\omega-\eta, \omega+\eta)$ we have $F(u)>\delta$,
that is to say
\[
 \pi(e^{u})-\li(e^{u}) > u^{-1}e^{u/2}\cdot\delta > 3.2 \times 10^{151}
\]
We have
\[
\li(N+r)-\li(N) - \frac{r}{\log N} = \int^{N+r}_{N}\frac{du}{\log
u} - \frac{r}{\log N} \leq 0 \leq \pi(N+r)-\pi(N)
\]
and so we obtain
\[
 \pi(N+r)-\li(N+r)\geq \pi(N)-\li(N) - \frac{r}{\log N}
\]
This shows that if $\pi(e^u)-\li(e^u) > M >0$ then $\pi(e^u) -
\li(e^u)$ will remain positive for another $[M u]$ consecutive
integers. Here, $Mu > 10^{154}$.

\begin{thm}There is a value of $u$ in the interval
\[(727.951858,727.952178)\] for which $\pi(e^{u}) - \li(e^{u}) >
3.2 \times 10^{151}$. There are at least $10^{154}$ successive
integers $x$ between $\exp(727.951858)$ and $\exp(727.952178)$ for
which $\pi(x)>\li(x)$.
\end{thm}

\end{document}